\documentclass[11pt]{amsart}
\usepackage[all]{xy}

\usepackage{amsmath}
\usepackage[latin1]{inputenc}
\usepackage{amsfonts}
\usepackage{amssymb}
\usepackage{multicol}
\usepackage{verbatim}
\usepackage{amsthm}
\usepackage{amscd}
\usepackage{pb-diagram}
\usepackage[mathscr]{eucal}
\usepackage[latin1]{inputenc}

\usepackage{amsmath}

\usepackage{amsfonts}

\usepackage{amssymb}

\usepackage{amsthm}

\usepackage{graphics}

\newcommand{\ZZ}{\mathbb{Z}}

\newcommand{\CC}{\mathbb{C}}

\newcommand{\Glie}{\mathfrak{g}}

\newcommand{\Yim}{\mathcal{Y}}

\newcommand{\Hlie}{\mathfrak{h}}

\newcommand{\U}{\mathcal{U}}

\newtheorem{thm}{Theorem}[section]

\newtheorem{rem}[thm]{Remark}

\newtheorem{ex}[thm]{Example}

\usepackage[all]{xy}
\title{
Quantum periodicity and Kirillov-Reshetikhin modules}
\author{David Hernandez}

\address{Universit\'e de Paris, Sorbonne Universit\'e, CNRS, Institut de
  Math\'ematiques de Jussieu-Paris Rive Gauche, F-75013 Paris, France}

\email{david.hernandez@imj-prg.fr}

\begin{document} 

\begin{abstract} We give a proof of the periodicity of quantum $T$-systems of type $A_n\times A_\ell$ with certain spiral boundary conditions.
Our proof is based on categorification of the $T$-system in terms of the representation theory of quantum affine algebras, more precisely on relations between classes of
Kirillov-Reshetikhin modules and of evaluation modules.
\end{abstract}

\maketitle

\begin{center}
{\it To Nicolai Reshetikhin on his 60th birthday}
\end{center}

\vskip 4.5mm

\tableofcontents

\section{Introduction}

The $Q$-system was introduced by Kirillov and Reshetikhin \cite{kr} as a system of relations between characters
of certain simple finite-dimensional representations of the quantum affine algebra $\mathcal{U}_q(\hat{sl}_{n+1})$, 
now called Kirillov-Reshetikhin modules 
$$Q_{a,b}^2 = Q_{a-1,b}Q_{a+1,b} + Q_{a,b+1}Q_{a,b-1},$$
where $1\leq a\leq n$ and $b$ is a non-negative integer.
Inspired by this work, the $T$-system was written in \cite{kns} as a refined version of the $Q$-system depending on 
a spectral parameter $u$ :
$$T_{a,b}(u-1)T_{a,b}(u+1) = T_{a-1,b}(u)T_{a+1,b}(u) + T_{a,b+1}(u)T_{a,b-1}(u).$$
It was conjectured that it is satisfied by the classes of the Kirillov-Reshetikhin modules.  

In a fundamental work \cite{Fre}, Frenkel and Reshetikhin introduced a character theory for finite-dimensional
representations of quantum affine algebras, called the $q$-characters. Then the $T$-system
is satisfied by $q$-characters of Kirillov-Reshetikhin modules in all types \cite{Nkr, Hkr}, and so by their classes as conjectured above.

In another direction, Zamolodchikov initiated in \cite{Z} a long series of work on the periodicity of solutions of $T$-systems 
with certain boundary conditions, which culminated in the work of Keller \cite{k} with a very general uniform proof of the 
periodicity of $T$-systems associated to a pair $(\Delta, \Delta')$ of Dynkin diagrams (see \cite{iikns, k1} for reviews and references).

In this note we propose a simple proof of the periodicity (and half-periodicity) of $T$-systems of type $A_n\times A_\ell$ and of its quantum version (in the sense of \cite{Nkr, HL:qGro}), 
with certain spiral boundary conditions (more general than the unit condition usually considered). We follow the approach in \cite[Section 12.1]{hl} where the proof of the commutative periodicity in type $A_n\times A_1$ is obtained with formulas for solutions in terms of $q$-characters. Indeed we find solutions in terms of 
certain evaluation representations, containing Kirillov-Reshetikhin modules but not only, and more precisely in terms of their $q,t$-characters defined by Nakajima \cite{Nak:quiver}.

The quantum periodicity (and half-periodicity) established in this note should also follow from the analog results 
in the commutative case (with unit boundary condition) mentioned above and from results in \cite{BZ, CKLP} 
(indeed the approach in \cite{k} is based on the study of the periodicity of a sequence of mutations in a certain cluster algebra). 
Our direct method gives an explicit solution in terms of $q,t$-characters.

The paper is organized as follows. In Section \ref{statement} we state the main periodicity and quantum periodicity results and we give several examples. In 
Section \ref{un} we give the necessary reminders on the representation theory of quantum affine algebras. In Section \ref{deux} we recall how the $T$-system
appears in the Grothendieck ring of the category of representations and we prove it has also other incarnations. We conclude in Section \ref{trois} with 
the proof of the quantum periodicity.

\medskip

{\bf Acknowledgment : } The author is very grateful to Bernard Leclerc for many discussions over the years and to Bernhard Keller
for useful remarks and explanations about \cite{k} and its consequences. The author would you like to thank Laura Fedele for her
careful reading and for pointing typos in a former version of this paper. The author is supported by the European Research Council under the European Union's Framework Programme H2020 with ERC Grant Agreement number 647353 Qaffine.

\section{Periodicity and quantum periodicity}\label{statement}

In this section we state the periodicity and quantum periodicity of $T$-systems that we establish in this note.

\subsection{Periodicity}

Let us first state the commutative $A_n\times A_\ell$-periodicity. Let 
$$I = \{1,\cdots, n\}\text{ and }J = \{1,\cdots, \ell\}.$$
We work on the lattice 
$$\Lambda = \{(a,b,u)\in I\times J\times \ZZ|a + b + u \in 2\ZZ\}.$$
Let us consider a family of commuting variables $(T_{a,b}(u))_{(a,b,u)\in \Lambda}$ satisfying the $T$-system 
(sometimes called octahedron relation) :
$$T_{a,b}(u-1)T_{a,b}(u+1) = T_{a-1,b}(u)T_{a+1,b}(u) + T_{a,b+1}(u)T_{a,b-1}(u),$$
for any $(a,b,u  + 1)\in \Lambda$.

So that the system is well-defined, we have to fix the boundary conditions, that is the values of 
$$T_{0,b}(u)\text{ , }T_{a,\ell + 1}(u)\text{ , }T_{n+1, b}(u)\text{ and }T_{a,0}(u).$$ 
The first choice is to set all values to $1$, this is called the unit boundary condition (see 
\cite{iikns}). The system is already non trivial with such a choice.
We will also consider the following boundary condition. Let $(\mathcal{F}_r)_{r\in I}$ be formal variables that we call coefficients.
We also set $\mathcal{F}_0 = \mathcal{F}_{n+1} = 1$. 
We set the following boundary conditions (for $u$ modulo $2(\ell + n + 2)$) : 
\begin{align*}
T_{n+1,m}(u) = \begin{cases} 
\mathcal{F}_{\frac{u + n  + m + 3}{2}} &\text{ for }    -n - 3 \leq u + m  \leq n -1
\\1 & \text{ for } 0 \leq u + m - n + 1\leq 2\ell + 2 .
\end{cases}
\end{align*}
\begin{align*}
T_{0,m}(u) = \begin{cases} \mathcal{F}_{\frac{u - m + 2}{2}} & \text{ for }  - 2\leq u - m \leq 2n ,
\\ 1 &\text{ for } 0 \leq u - m - 2n   \leq 2 (\ell + 1)  .\end{cases}
\end{align*}
\begin{align*}
T_{k,\ell + 1}(u) = \begin{cases}\mathcal{F}_{\frac{u - k - \ell + 1}{2}}&\text{ for } 0\leq u - \ell  - k   + 1\leq 2n  + 2,
\\ 1
&\text{ for }  \ell +  1 \leq u - 2n - k \leq  3 (\ell + 1) .\end{cases}
\end{align*}
\begin{align*}
T_{k,0}(u) = \begin{cases} \mathcal{F}_{\frac{u + k + 2}{2}}& \text{ for }-2\leq u + k   \leq  2 n ,
\\ 1 &\text{ for } 0 \leq u  - 2n + k \leq 2(\ell   + 1  )  .\end{cases}
\end{align*}
This is a particular case of the 
spiral boundary condition (see \cite{iikns}).

\begin{rem}\label{initc} For $(a,b,u)\in \Lambda$, an induction on $u\geq a + b - 2$ shows that $T_{a,b}(u)$ is a rational fraction in the 
$$X_{k,m} = T_{k,m}(k + m - 2)\text{ for $(k,m)\in I\times J$}$$
and in the coefficients 
$$X_{k,0} = T_{k,0}(k - 2) = \mathcal{F}_k.$$
\end{rem}

We have the following periodicity.

\begin{thm}\label{clasper} For any $(a,b,u)\in \Lambda$, we have the half-periodicity property : 
$$T_{a,b}(u) = T_{n + 1 - a,\ell + 1 - b}(u + n + \ell + 2).$$
It implies that $T_{a,b}(u)$ is $2(n + \ell + 2)$-periodic in $u$.
\end{thm}

\begin{ex}\label{11} Let $n = \ell = 1$. The non trivial boundary conditions are
$$T_{0,1}(1) = T_{1,2}(3) = T_{2,1}(5) = T_{1,0}(7) = \mathcal{F}_1.$$
Set $X = X_{1,1} = T_{1,1}(0)$. 
\\The values of $T_{1,1}(u)$ for $u = 0, 2, 4, 6 , 8$ are respectively
$$X\text{ , }\frac{\mathcal{F}_1 + 1}{X}\text{ , }X\text{ , }\frac{\mathcal{F}_1 + 1}{X}\text{ , }X.$$
\end{ex}

\begin{ex}\label{12} Let $n = 1$ and $\ell = 2$. The non trivial boundary conditions are
$$T_{0,1}(1) = T_{0,2}(2) = T_{1,3}(4) = T_{2,2}(6) = T_{2,1}(7) = T_{1,0}(9) = \mathcal{F}_1.$$
Set $X_1 = X_{1,1} = T_{1,1}(0)$ and $X_2 = X_{1,2} = T_{1,2}(1)$. 
\\The values of $T_{1,1}(t)$ for $t = 0, 2, 4, 6 , 8, 10$ are
$$X_1\text{ , }\frac{\mathcal{F}_1 + X_2}{X_1}\text{ , }\frac{X_1 + 1 }{X_2}\text{ , }X_2\text{ , }\frac{\mathcal{F}_1 X_1 + \mathcal{F}_1 + X_2}{X_1X_2}\text{ , }X_1.$$
The values of $T_{1,2}(u)$ for $u = 1, 3, 5, 7 , 9, 11$ are respectively
$$X_2\text{ , }\frac{\mathcal{F}_1 X_1 + \mathcal{F}_1 + X_2}{X_1X_2}\text{ , }X_1\text{ , }\frac{\mathcal{F}_1 + X_2}{X_1}\text{ , }\frac{X_1 + 1 }{X_2}\text{ , }X_2.$$
\end{ex}

\begin{ex}\label{21} Let $n = 2$ and $\ell = 1$. The non trivial boundary conditions are
$$T_{0,1}(1) = T_{1,2}(3)  = T_{2,2}(4) = T_{3,1}(6) = T_{2,0}(8) = T_{1,0} (9) =   \mathcal{F}_1 $$
$$T_{1,0}(1) = T_{0,1}(3)  =  T_{1,2}(5) = T_{2,2}(6) = T_{3,1}(8) = T_{2,0}(11) = \mathcal{F}_2$$
Set $X_1 = X_{1,1} = T_{1,1}(0)$ and $X_2 = X_{2,1} = T_{2,1}(1)$. 
\\The values of $T_{1,1}(t)$ for $t = 0, 2, 4, 6 , 8, 10$ are
$$X_1\text{ , }\frac{\mathcal{F}_1 X_2 + \mathcal{F}_2}{X_1}\text{ , }\frac{X_1 + \mathcal{F}_2 }{X_2}\text{ , }X_2
\text{ , }\frac{X_1 + \mathcal{F}_2 + \mathcal{F}_1 X_2}{X_1X_2}\text{ , }X_1.$$
The values of $T_{2,1}(u)$ for $u = 1, 3, 5, 7 , 9, 11$ are respectively
$$X_2\text{ , }\frac{X_1 + \mathcal{F}_2 + \mathcal{F}_1X_2}{X_1X_2}\text{ , }X_1\text{ , }\frac{\mathcal{F}_2 + \mathcal{F}_1 X_2}{X_1}\text{ , }\frac{X_1 + \mathcal{F}_2 }{X_2}\text{ , }X_2.$$
\end{ex}

\subsection{Quantum periodicity}\label{qps}

Let us now the state the quantum version of the $A_n\times A_\ell$-periodicity (see also \cite{kn} and \cite{dfk} for $n = 1$).

We work now with quasi-commuting variables $(X_{a,b})_{(a,b)\in I\times J}$ :  
$$X_{a,b} * X_{c,d} = t^{\gamma(a,b;c,d) - \gamma(c,d; a,b)}X_{c,d} * X_{a,b}.$$
To define the power of $t$, we use the inverse $\tilde{C}(z)$ of the quantized Cartan matrix 
$$C(z) = ((z + z^{-1})\delta_{i,j} - \delta_{i+1,j} - \delta_{i-1,j})_{i,j\in I}.$$
For $p\in\mathbb{Z}$ and $a,c\in I$, we denote by $\tilde{C}_{a,c}(p)$ the coefficient of $z^p$ in the expansion in $z$ of $\tilde{C}_{a,c}(p)$. We set
$$\gamma(a,b;c,d) = \tilde{C}_{a,c}(2\ell - 2 b + c - a +1) +  \tilde{C}_{a,c}(2\ell - 2b + c - a  - 1)+ \cdots + \tilde{C}_{a,c}(2d - 2b + c - a + 1).$$

The relation is also extended to $b = 0$ or $d = 0$ so that we get the quasi-commutation rule with the coefficients. 
In particular for $r < r'$, we have 
\begin{equation}\label{ft}\mathcal{F}_r * \mathcal{F}_{r'} = t^{\tilde{C}_{r,r'}(2\ell + r' - r +1) +\tilde{C}_{r,r'}(2\ell + r' - r - 1)+\cdots + \tilde{C}_{r,r'}(2\ell + r - r' +3)}\mathcal{F}_{r'}*\mathcal{F}_r\end{equation}
and $\mathcal{F}_r * \mathcal{F}_{r+1} = t^{\tilde{C}_{r,r+1}(2\ell + 2)} \mathcal{F}_{r+1}* \mathcal{F}_r$. 
This is derived from $\tilde{C}_{i,j}(k) = 0$ if $k\leq |j - i|$ (which can be observed for example in the formula in \cite[Appendix A.3]{gtl}). 

The quasi-commuting variables $(X_{a,b})_{(a,b)\in I\times J}$ with the $\mathcal{F}_t$ generate a quantum torus $\mathcal{T}_t$ over $\mathbb{Z}[t^{\pm 1/2}]$. 
We denote its fraction field by $K_t$. It has an antimultiplicative bar involution satisfying $\overline{t} = t^{-1}$ and so that the $X_{a,b}$, $\mathcal{F}_t$ are bar-invariant.

For each product $m$ of various $X_{a,b}^{\pm 1}$, $\mathcal{F}_t^{\pm 1}$, there is a unique $\alpha\in\ZZ$ so that $t^{\alpha/2}m$ is bar invariant. This is called a commutative
monomial. The commutative monomials form a $\mathbb{Z}[t^{\pm 1/2}]$-basis of $\mathcal{T}_t$. 

\begin{ex} For $n = \ell = 1$, we have : 
$$X_1*\mathcal{F}_1 = t^2 \mathcal{F}_1 * X_1.$$
For $n = 1$, $\ell = 2$, we have : 
$$X_1 * X_2 = t^{-2} X_2 * X_1\text{ , }X_1*\mathcal{F}_1 = \mathcal{F}_1 * X_1\text{ , }X_2 * \mathcal{F}_1 = \mathcal{F}_1 * X_2.$$
For $n = 2$, $\ell = 1$, we have : 
$$X_1 * X_2 = t X_2 * X_1\text{ , }X_1 * \mathcal{F}_1 = t \mathcal{F}_1 * X_1\text{ , }X_1 * \mathcal{F}_2 = \mathcal{F}_2 * X_1,$$
$$X_2 * \mathcal{F}_2 = t \mathcal{F}_2 * X_2\text{ , }X_2 * \mathcal{F}_1 = t \mathcal{F}_1 * X_2\text{ , }\mathcal{F}_1 * \mathcal{F}_2 = t^{-1} \mathcal{F}_2*\mathcal{F}_1.$$
\end{ex}

We fix the same the boundary conditions as for the commutative setting above.

\begin{thm}\label{qperiodic} Consider a family of bar-invariant $T_{a,b}(u)\in K_t$ satisfying : 
$$T_{a,b}(u - 1) * T_{a,b}(u+1) \in t^{\ZZ/2}  T_{a-1,b}(u) * T_{a+1,b}(u) + t^{\ZZ/2} T_{a,b+1}(u) * T_{a,b-1}(u)$$
for $(a,b,u+1)\in \Lambda$. We assume the same initial conditions as in Remark \ref{initc} and the same spiral boudary conditions as in the classical setting. Then for $(a,b,u)\in \Lambda$ :
$$T_{a,b}(u) = T_{n + 1 - a,\ell + 1 - b}(u + n + \ell + 2).$$
It implies that $T_{a,b}(u)$ is $2(n + \ell + 2)$-periodic in $u$.
\end{thm}

The classical periodicity in Theorem \ref{clasper} follows directly from this Theorem.
We propose a simple proof based on the representations theory of quantum affine algebras and 
on their Kirillov-Reshetikhin modules.

\begin{ex}
Let us study the examples \ref{11}, \ref{12}, \ref{21} above. In these examples, let us just replace each 
Laurent monomial in the $X_{a,b}$ by the corresponding commutative monomial in the quantum torus. We get 
bar-invariant elements in $\mathcal{T}_t$ and we keep the notation $T_{a,b}(u)$. Let us verify they satisfy the quantum $T$-system.

Let $n = \ell = 1$. We get : 
\begin{align*}
T_{1,1}(0) * T_{1,1}(2) = t \mathcal{F}_1 + 1&, \quad T_{1,1}(2) * T_{1,1}(4) = t^{-1}\mathcal{F}_1 + 1.\end{align*}
Let $n = 1$, $\ell = 2$. 
\begin{align*}
T_{1,1}(0) * T_{1,1}(2) = \mathcal{F}_1 + t^{-1} T_{1,2}(1)&, \quad T_{1,2}(1) * T_{1,2}(3) = \mathcal{F}_1 + t^{-1}T_{1,1}(2),
\\T_{1,1}(2) * T_{1,1}(4) = 1 + t^{-1}T_{1,2}(3)&, \quad  T_{1,2}(3) * T_{1,2}(5) = 1 + t^{-1} \mathcal{F}_1 * T_{1,1}(4),
\\T_{1,1}(4) * T_{1,1}(6) = 1 + t^{-1}T_{1,2}(5)&, \quad T_{1,2}(5) * T_{1,2}(7) = \mathcal{F}_1 + t^{-1} T_{1,1}(6),
\\T_{1,1}(6) * T_{1,1}(8) = \mathcal{F}_1 + t^{-1}T_{1,2}(7)&, \quad T_{1,2}(7) * T_{1,2}(9) = 1 + t^{-1}T_{1,1}(8),
\\T_{1,1}(8) * T_{1,1}(10) = 1 + t^{-1}T_{1,2}(9)*\mathcal{F}_1&, \quad T_{1,2}(9) * T_{1,2}(11) = 1 + t^{-1}T_{1,1}(10).
\end{align*}
Let $n = 2$, $\ell = 1$. 
\begin{align*}
T_{1,1}(0) * T_{1,1}(2) = t^{\frac{1}{2}}T_{2,1}(1) * \mathcal{F}_1 + \mathcal{F}_2&, \quad 
T_{2,1}(1) * T_{2,1}(3) = t T_{1,1}(2) + 1,
\\T_{1,1}(2) * T_{1,1}(4) = t^{\frac{3}{2}} T_{2,1}(3) * \mathcal{F}_2 +  \mathcal{F}_1&, \quad 
T_{2,1}(3) * T_{2,1}(5) = t^{\frac{1}{2}} T_{1,1}(4) + t^{-\frac{1}{2}}\mathcal{F}_1,
\\T_{1,1}(4) * T_{1,1}(6) = t^{\frac{1}{2}}T_{2,1}(5) + t^{-\frac{1}{2}}\mathcal{F}_2&, \quad 
T_{2,1}(5) * T_{2,1}(7) = t^{\frac{3}{2}}\mathcal{F}_1 * T_{1,1}(6)  + \mathcal{F}_2,
\\T_{1,1}(6) * T_{1,1}(8) = t T_{2,1}(7) + 1&, \quad 
T_{2,1}(7) * T_{2,1}(9) = t^{\frac{1}{2}} \mathcal{F}_2 * T_{1,1}(8)  +  \mathcal{F}_1,
\\T_{1,1}(8) * T_{1,1}(10) = t^{\frac{1}{2}} T_{2,1}(9) + t^{-\frac{1}{2}}\mathcal{F}_1&, \quad 
T_{2,1}(9) * T_{2,1}(11) = t^{\frac{1}{2}}T_{1,1}(10) + t^{-\frac{1}{2}}\mathcal{F}_2.
\end{align*}
\end{ex}

\section{Finite-dimensional representations of quantum affine algebras}\label{un}

We recall the main definitions and properties of finite-dimensional representations of the quantum affine algebra associated to $sl_{n+1}$. 

\subsection{Quantum affine algebras}  

All vector spaces, algebras and tensor products are defined over $\CC$.

Let $C = (C_{i,j})_{0\leq i,j\leq n}$ be the Cartan matrix of type $A_n^{(1)}$, that is 
$$C_{i,j} = 2\delta_{i,j} - \delta_{i,j+1} - \delta_{i+1,j}$$
where $n+1$ is identified with $0$. Fix $q\in\CC^*$ which is not a root of unity.

The {\it quantum affine algebra} $\U_q(\Glie)$ is defined by generators $k_i^{\pm 1}$, $x_i^{\pm}$ ($0\leq i\leq n$) and 
relations
$$k_ik_j=k_jk_i\text{ , } k_ix_j^{\pm}=q^{\pm C_{i,j}}x_j^{\pm}k_i\text{ , }[x_i^+,x_j^-]=\delta_{i,j}\frac{k_i-k_i^{-1}}{q-q^{-1}},$$
$$\underset{p=0\cdots 1 - C_{i,j}}{\sum}(-1)^p(x_i^{\pm})^{\left(1-C_{i,j}-p\right)}x_j^{\pm}(x_i^{\pm})^{(p)}=0 \text{ (for $i\neq j$)}
,$$
where we denote $\left(x_i^{\pm}\right)^{(p)} = \left(x_i^{\pm}\right)^p/[p]_{q}$ for $0\leq p\leq 2$, where $[p]_q = (q^p - q^{-p}) (q - q^{-1})^{-1}$. 

It is a Hopf algebra with a coproduct $\Delta : \U_q(\Glie)\rightarrow \U_q(\Glie)\otimes \U_q(\Glie)$ defined for $0\leq i\leq n$ by
$$\Delta(k_i)=k_i\otimes k_i\text{ , }\Delta(x_i^+)=x_i^+\otimes 1 + k_i\otimes x_i^+\text{ , }\Delta(x_i^-)=x_i^-\otimes k_i^{-1} + 1\otimes x_i^-.$$ 

Let $\overline{\Glie} = sl_{n+1}$ be the finite-dimensional simple Lie algebra of Cartan matrix $(C_{i,j})_{i,j\in I}$. We denote respectively by $\omega_i$, $\alpha_i$, $\alpha_i^\vee$ ($i\in I$) the fundamental weights, the simple roots and the simple coroots of $\overline{\Glie}$.
We use the standard partial ordering $\leq$ on the weight lattice $P$ of $\overline{\Glie}$. 

The algebra $\U_q(\Glie)$ has another set of generators, the {\it Drinfeld generators}, denoted by 
$$x_{i,m}^\pm\text{ , }k_i^{\pm 1}\text{ , }h_{i,r}\text{ , }c^{\pm 1/2}\text{ for $i\in I$, $m\in \ZZ$, $r\in\ZZ \setminus \{0\}$.}$$ 
We have $x_i^\pm = x_{i,0}^\pm$ for $i\in I$. A complete set of relations for Drinfeld generators was obtained in \cite{bec, da2}. In particular the multiplication defines a surjective linear morphism
\begin{equation}\label{trian}\U_q^-(\Glie)\otimes \U_q(\Hlie)\otimes \U_q^+(\Glie)\rightarrow \U_q(\Glie)\end{equation}
where $\U_q^\pm(\Glie)$ is the subalgebra generated by the $x_{i,m}^\pm$ ($i\in I$, $m\in\ZZ$) and $\U_q(\Hlie)$ is the subalgebra generated by the $k_i^{\pm 1}$, the $h_{i,r}$ and $c^{\pm 1/2}$ ($i\in I$, $r\in\ZZ\setminus\{0\}$).

\subsection{Finite-dimensional representations}\label{fdrap}
We refer to \cite{CH} for generalities on the category $\mathcal{C}$ of finite-dimensional representations of $\U_q(\Glie)$.
For $i\in I$, the action of $k_i$ on any object of $\mathcal{C}$ is diagonalizable with eigenvalues in $\pm q^{\ZZ}$.
Without loss of generality, we can assume that $\mathcal{C}$ is the category of {\it type 1} finite-dimensional representations (see \cite{Cha2}), i.e. we assume that for any object of $\mathcal{C}$, the eigenvalues of $k_i$ are in $q^{\ZZ}$ for $i\in I$. The simple objects of $\mathcal{C}$ are parametrized by $n$-tuples of polynomials $(P_i(u))_{i\in I}$ satisfying $P_i(0) = 1$ (they are called {\it Drinfeld polynomials}) \cite{cp, Cha2}. 

In type $A$, there is a family of evaluation morphisms $ev_a : \U_q(\Glie)\rightarrow \U_q(\overline{\Glie})$ parametrized by $a\in\CC^*$. 
Hence for $V$ a simple finite-dimensional representations of $\U_q(\overline{\Glie})$, by pull-back we get an evaluation representation $(V)_a$.
If the highest weight of $V$ is a multiple of a fundamental weight, then $V$ is a Kirillov-Reshetikhin module. In the particular case of a fundamental weight,
we get the fundamental representations $V_i(a) = (V(\omega_i))_a$ of Drinfeld polynomials $(1,\cdots, 1, 1 - za , 1, \cdots, 1)$ with a non-trivial polynomial in position $i$. Their classes generate the Grothendieck ring $K_0(\mathcal{C})$ of the category $\mathcal{C}$ which
is a polynomial ring in the variables $[V_i(a)]$ as proved in \cite{Fre}. In general simple finite-dimensional representations are not evaluation modules.

For $\omega\in P$, the {\it weight space} $V_\omega$ of an object $V$ in $\mathcal{C}$ is the set of {\it weight vectors} of weight $\omega$, i.e. of vectors $v\in V$ satisfying $k_i v = q^{\left( \omega(\alpha_i^\vee)\right)} v$ for any $i\in I$. 

The elements $c^{\pm 1/2}$ act by identity on any object $V$ of $\mathcal{C}$, and so the action of
the $h_{i,r}$ commute. Since the $h_{i,r}$, $i\in I$, $r\in\mathbb{Z}\setminus \{0\}$, also commute with the $k_i$, $i\in I$, every object in $\mathcal{C}$ can be decomposed as a direct sum of generalized eigenspaces of the $h_{i,r}$ and $k_i$. More precisely, by Frenkel-Reshetikhin theory of $q$-characters \cite{Fre}, the eigenvalues of the $h_{i,r}$ and $k_i$ can be {\it encoded} by {\it monomials} $m$ in formal variables $Y_{i,a}^{\pm 1}$ ($i\in I, a\in\CC^*$). Let $\mathcal{M}$ be the set of such monomials (also called {\it $l$-weights}). Given $m\in\mathcal{M}$ and an object $V$ in $\mathcal{C}$, 
let $V_m$ be the subspace of $V$ of 
common pseudo-eigenvectors of the $h_{i,r}$, $k_i$ with pseudo-eigenvalues associated to $m$ 
(also called {\it $l$-weight space}). Thus, 
$$V = \bigoplus_{m\in\mathcal{M}} V_m.$$  
If $v\in V_m$, then $v$ is a weight vector of weight 
$$\omega(m) = \sum_{i\in I, a\in\CC^*} u_{i,a}(m) \omega_i\in P,$$ 
where we denote $m = \prod_{i\in I, a\in\CC^*}Y_{i,a}^{u_{i,a}(m)}$. For $v\in V_m$, we set $\omega(v) = \omega(m)$.

The {\it $q$-character morphism} is an injective ring morphism
$$\chi_q : \text{Rep}(\U_q(\Glie)) \rightarrow \Yim = \ZZ\left[Y_{i,a}^{\pm 1}\right]_{i\in I, a\in\CC^*},$$
$$\chi_q(V) = \sum_{m\in \mathcal{M}} \text{dim}(V_m) m.$$
If $V_m\neq \{0\}$ we say that $m$ is an {\it $l$-weight of $V$}.

A monomial $m\in\mathcal{M}$ is said to be {\it dominant} if $u_{i,a}(m)\geq 0$ for any $i\in I, a\in\CC^*$. For $V$ a simple object in $\mathcal{C}$, let $M(V)$ be the {\it highest weight monomial} of $\chi_q(V)$, that is so that $\omega(M(V))$ is maximal for the partial ordering on $P$. $M(V)$ is dominant and characterizes the isomorphism class of $V$ (it is equivalent to the data of the Drinfeld polynomials). Hence to a dominant monomial $M$ is associated a simple representation $L(M)$. For $i\in I$ and $a\in\CC^*$, we have for example the fundamental representation $V_i(a) = L\left(Y_{i,a}\right)$. The simple modules of highest weight monomial
$$X_{i,\alpha}^\beta= Y_{i,q^{\alpha}} Y_{i,q^{\alpha + 2}}\cdots Y_{i,q^{\alpha + 2 (\beta - 1)}}$$
for some $i\in I, \alpha\in\ZZ, \beta\geq 1$ are Kirillov-Reshetikhin modules. We will also use the notation $X_{i,\alpha}^\beta = 1$ for $\beta \leq 0$.

\begin{ex} The $q$-character of the fundamental representation $L(Y_a)$ of $\U_q(\hat{sl_2})$ is 
$$\chi_q(L(Y_a)) = Y_a + Y_{aq^2}^{-1}.$$
\end{ex}

The $q$-characters of evaluation modules, including Kirillov-Reshetikhin modules and fundamental modules, are known explicitly 
(see references in the introduction of \cite{miniaff}). The formulas involve
the monomials $A_{i,a}$ defined in \cite{Fre} for $i\in I, a\in \CC^*$ by
$$A_{i,a} = Y_{i,aq^{-r_i}}Y_{i,aq^{r_i}}\times \prod_{\left\{j\in I|C_{i,j} = -1\right\}}Y_{j,a}^{-1}.$$

\subsection{Quantum Grothendieck ring} The Grothendieck ring $K_0(\mathcal{C})$ has a $t$-deformation called quantum Grothendieck ring 
$K_t(\mathcal{C})$ as constructed in \cite{VV:qGro, Nak:quiver, H1} (we use the version of \cite{H1, HL:qGro}).
It is a $\ZZ[t^{\pm 1/2}]$-subalgebra of a quantum torus $\mathcal{Y}_t$ and simple objects $L(m)$ have corresponding classes $[L(m)]_t\in K_t(\mathcal{C})$.
A quantum version of a result in \cite{Fre2} gives the following \cite{Nak:quiver} :
\begin{equation}\label{order}[L(m)]_t \in m * ( 1 +  \ZZ[t^{\pm 1/2}, A_{i,c}^{-1}]_{i\in I, c\in\CC^*}).\end{equation}
In other words, $m$ is maximal for the Nakajima partial ordering on monomials, that is $M\preceq M'$ if $M'M^{-1}$ is a product of variables $A_{i,c}$.

If a simple module $V$ is thin, that is if its $\ell$-weight spaces are of dimension $1$, then $[V]_t$ is a sum of commutative
monomials (defined as in section \ref{qps}) and can be identified with its $q$-character (see \cite[Corollary 5.3]{HL:qGro}). In type $A$, all simple evaluation
modules are thin.


\section{Relations in the Grothendieck ring}\label{deux}

We recall how the $T$-system originally occurs in the representation theory of quantum affine algebras and we also
establish another incarnation of the $T$-system in the Grothendieck ring $K_0(\mathcal{C})$ (that we call horizontal $T$-system).

We denote $I = \{1,\cdots, n\}$ and $J = \{1,\cdots, \ell\}$ as above.

\subsection{Original $T$-systems}
For $1\leq i\leq n$ and $0\leq m \leq p\leq \ell$, consider the Kirillov-Reshetikhin module
$$\beta(m,p)^i = L(X_{i,i+2m}^{p-m+1}).$$
We extend the notation to $i = 0$ and $i = n + 1$ by setting $\beta(m,p)^i  = 1$ in these cases.

For $i\in I$ and $0\leq m \leq p < \ell$, we have the $T$-system in $K_0(\mathcal{C})$:
\begin{equation*}\label{betal}\beta(m,p)^i\beta(m+1,p+1)^i = \beta(m,p+1)^i\beta(m+1,p)^i + \beta(m+1,p+1)^{i-1}\beta(m,p)^{i+1}.\end{equation*}
See the list of references in the introduction of \cite{miniaff}. It can be deformed into the quantum $T$-system in $K_t(\mathcal{C})$ (see \cite{Nkr, HL:qGro}) :
\begin{equation}\label{tbetal}\beta(m,p)^i*\beta(m+1,p+1)^i = t^\lambda \beta(m,p+1)^i*\beta(m+1,p)^i + t^\mu \beta(m+1,p+1)^{i-1}*\beta(m,p)^{i+1}.\end{equation}
for some $\lambda,\mu\in\mathbb{Z}/2$ which depend in $m,p, i$ (they can be explicitly computed but this is not relevant for the following).

\subsection{Horizontal $T$-systems}

The $T$-system has another incarnation in $K_0(\mathcal{C})$. 

For $0\leq i \leq j\leq n+1$ and $0\leq m\leq \ell + 1$, consider the evaluation module
$$\alpha(i,j)^m = L(M_{[i,j]}^m)\text{ where }M_{[i,j]}^m = X_{i,i}^{m}X_{j,j+2m}^{\ell + 1 - m}.$$
Some of these representations are Kirillov-Reshetikhin modules : 
\begin{equation}\label{exkr}\alpha(i,n+1)^{m+1} = \beta(0,m)^i\text{ and }\alpha (0,i)^{m} = \beta(m,\ell)^i\text{ for $0\leq m\leq \ell$ 
and $0\leq i\leq n+1$.}\end{equation}
For $0\leq i\leq n+1$ and $j\geq 0$, we will denote
\begin{equation}\label{fi}F_i = L(X_{i,i}^{\ell + 1}) = \alpha(i,i)^m = \beta(0,\ell)^i = \alpha(i,i + j)^{\ell + 1} = \alpha(i-j,i)^0.\end{equation}
Note that $F_0 = F_{n+1} = 1$.
\begin{thm}\label{typeahl} For $0\leq i < j \leq n$ and $m\in J$, there are $\lambda, \lambda '\in\ZZ/2$, so that :
$$\alpha(i,j)^m* \alpha(i+1,j+1)^m
= t^\lambda \alpha(i,j+1)^m*\alpha(i+1,j)^m + t^{\lambda'}\alpha(i,j)^{m+1} * \alpha(i+1,j+1)^{m-1}.$$
\end{thm}

\begin{rem} (i) This relation is "orthogonal" to the original quantum $T$-system in the sense that the spectral parameter 
is replaced by the vertex of the Dynkin diagram.

(ii) At $t = 1$, it can be shown that the relation comes from a non-split exact sequence obtained by a normalized $R$-matrix, as for the
original $T$-system (see \cite{Nkr, Hkr}). In fact, it can be checked that the two tensor products associated to the right hand terms correspond to simple modules.
Using \cite[Theorem 4]{c}, the proof is analogous to the one for the $T$-system.

(iii) At $t = 1$, this relation can be seen as an extended $T$-systems in \cite{my2}. 

(iv) For the limit values of $i,j$, the relation involves both the $\alpha$ and the $\beta$-modules and so connect the two families. 
The specialization at $t = 1$ reads : 
$$\beta(m,\ell)^j\alpha(1,j+1)^m
= \beta(m,\ell)^{j+1}\alpha(1,j)^m + \beta(m+1,\ell)^j \alpha(1,j+1)^{m-1}\text{( for $i = 0$)}$$
$$\alpha(i,n)^m\beta(0,m-1)^{i+1}
= \beta(0,m-1)^i\alpha(i+1,n)^m + \alpha(i,n)^{m+1} \beta(0,m-2)^{i+1}\text{( for $j = 0$)}.$$
\end{rem}

\begin{proof} By \cite{Fre, Fre2}, a $q$-character is determined uniquely by the multiplicity of its dominant monomials. 
We will use the notation $A_{i,\lambda}$ instead of $A_{i,q^\lambda}$ for $i\in I$, $\lambda\in\mathbb{Z}$.
First we prove that $\alpha(i,j)^m* \alpha(i+1,j+1)^m$ has $2m+1$ dominant monomials :
$$M_1 = M_{[i,j]}^m M_{[i+1,j+1]}^m\text{ , }M_2 = M_1 A_{i+1,i+2m}^{-1}A_{i+2,i+1+2m}^{-1}
\cdots A_{j,j+2m-1}^{-1},$$
$$M_{2r} = M_2 \prod_{2\leq p\leq r} (A_{i,i+2m-2p+3}A_{i+1,i+2m-2p+2})^{-1}\text{ , }M_{2r+1} = M_{2r} A_{i,i+2m-2r+1}^{-1},$$ 
where $1\leq r \leq m$. It is clear that these monomials occur. Indeed $M_1$ is the product of the highest monomials. 
For $2\leq r\leq 2m+1$, we decompose 
$$M_r = (M_{[i+1,j+1]}^m (M_2M_1^{-1}))\times (M_{[i,j]}^m (M_rM_2^{-1})).$$
Now consider $M\neq M_1$ a dominant monomial which occurs. We factorize 
$$M = M_1 M' M''$$
where $M'$ (resp. $M''$) is a monomial of $(M_{[i,j]}^m)^{-1}\alpha(i,j)^m$ 
(resp. of $(M_{[i+1,j+1]}^m)^{-1}\alpha(i+1,j+1)^m$). 
As $M$ is dominant, we have $M'M''\in\ZZ[A_{k,r}^{-1}]_{k\in I, r\leq j + 2\ell}$. Then
$$M'\in \ZZ[A_{k,r}^{-1}]_{k\leq j-1,r\in\ZZ}$$
and that there is $R\geq 0$ such that
$$M''\in (A_{j,j+2m-1}A_{j,j+2m-3} \cdots A_{j,j+2m-1-2R})^{-1} \ZZ[A_{k,r}^{-1}]_{k\leq j-1,r\in\ZZ}.$$
The monomial $\tilde{M} = M (X_{j+1,j+1+2m}^{\ell + 1 - m}X_{j,j+2m}^{\ell + 1 - m})^{-1}$ 
is a monomial of $\chi_q(L(X_{i,i}^mX_{i+1,i+1}^m))$ which has a unique dominant monomial. 
Hence $\tilde{M} Y_{j,j+2m}$ is dominant. So 
$$\tilde{M}' = \tilde{M}(A_{j,j+2m-1}\cdots A_{i+1,i+2m}^{-1})$$ 
is a monomial of 
$\chi_q(L(X_{i,i}^mX_{i+1,i+1}^m))$. If $\tilde{M}' = X_{i,i}^mX_{i+1,i+1}^m$,
then $M = M_2$. Otherwise, $\tilde{M}A_{i,i+2m-1}$ is a monomial
of $\chi_q(L(X_{i,i}^mX_{i+1,i+1}^m))$. We continue by induction, and so 
$M$ is one of the $M_r$.

This also implies that each $M_r$ occurs with multiplicity which is a power of $t$. 

Similarly, we get that $\alpha(i,j+1)^m * \alpha(i+1,j)^m$ has $m+1$
dominant monomials which are the $M_{2r+1}$ for $0\leq r\leq m$.
We also get that $\alpha(i,j)^{m+1} * \alpha(i+1,j+1)^{m-1}$ has $m$
dominant monomials which are the $M_{2r}$ for $1\leq r\leq m$.

To conclude, we have to check that the powers of $t$ match : this can be done using positivity in the quantum Grothendieck ring as in \cite[Section 5.10]{HL:qGro} 
or directly as in \cite[section 9]{ho}.
\end{proof}

\section{Proof of periodicity}\label{trois} 

In this section we finish the proof of the quantum periodicity.

It suffices to identify the $T_{a,b}(t)$ with variables
satisfying the $T$-system, the half-periodicity and such that the variables corresponding to the $X_{a,b}$ are
algebraically independent.
We will identify the $T_{a,b}(t)$ with certain $q,t$-characters of minimal affinizations, that is elements of the quantum torus $\mathcal{Y}_t$.

For $0\leq k\leq n+1$, $0\leq m\leq \ell+1$ and $u\in\ZZ$ so that $k + m + u\in 2\ZZ$, we set :
\begin{align*}
T_{k,m}(u) = \begin{cases} \alpha(\frac{u + 2 - k - m}{2},\frac{u + 2 + k - m}{2})^m & \text{ for }0\leq u +2 - k - m \leq 2(n + 1 - k),
\\ \beta(\frac{u - 2n + k - m}{2},\frac{u - 2n - 2 + k + m}{2})^{n + 1 - k}&\text{ for } m \leq u -2n + k   \leq 2\ell - m  + 2 ,
\\ \alpha(\frac{u - 2n - 2\ell + k + m - 2}{2}, \frac{u - k - 2 \ell + m }{2})^{\ell + 1 - m}&\text{ for } 0\leq u  -2n -2\ell - 2 + m +  k \leq 2k ,
\\ \beta(\frac{u - 2 - 2n - k + m - 2\ell}{2},\frac{u - 2n - 2  - k - m}{2})^k
&\text{ for }- m \leq u - 2\ell   -2n - k - 2  \leq  m  .\end{cases}
\end{align*}
This defines $T_{k,m}(u)$ for $0\leq u - m - k + 2\leq 2n + 2\ell + 4$, and we extend the definition for any $u$ by $2(n + \ell + 2)$-periodicity.

\begin{rem}
(i) The formulas in all cases are compatible thanks to relations (\ref{exkr}).

(ii) Identifying the class $F_r$ defined in (\ref{fi}) with $\mathcal{F}_r$, we recover boundary conditions of Section \ref{statement}.
\end{rem}

The $X_{k,m}$ quasi-commute, with the same rules as in Section \ref{statement}.
The relations (\ref{tbetal}) and Theorem \ref{typeahl} imply that the $T_{a,b}(u)$ satisfy the quantum $T$-system
for a distinguished choice of the powers of $t$ (let us call it the distinguished powers). 

The $(X_{k,m})_{(k,m)\in I\times (J\cup\{0\})}$ form a family of 
algebraically independent variables. We may argue as in \cite{hl2}. Let us explain this point for completeness : all the representations we consider belong to the 
monoidal category $\mathcal{C}_\ell^o$ of representations whose classes belong to the subring of the Grothendieck ring $K_0(\mathcal{C})$
generated by the classes of fundamental representations $[L(Y_{k,k+2m})]$ for $(k,m)\in I\times (J\cup\{0\})$.
Then there is an injective ring morphism 
$$\chi_q^T : K_0(\mathcal{C}_\ell^o) \rightarrow \mathcal{Y}$$
called truncated $q$-character morphism \cite{hlad} : it is defined so that for $L(m)$ a simple module in $\mathcal{C}_\ell^o$, 
$\chi_q^T(L(m))$ is obtained from $\chi_q(L(m))$ by removing the monomials $m'$ so that in $m'm^{-1}$ contains a factor 
of the form $A_{k,k+2\ell + 1}^{-1}$, $k\in I$. Now by \cite{Hkr}, the $\chi_q^T(X_{k,m}) = X_{k,k+2m}^{\ell + 1 - m}$ are just monomials which are clearly
algebraically independent.

As by construction we have $T_{a,b}(u) = T_{n + 1 - a,\ell + 1 - b}(u + n + \ell + 2)$, we get the result for 
the quantum $T$-system with the distinguished powers of $t$.

To conclude it suffices to check that the powers of $t$ correspond automatically to the distinguished choice. 
We consider a solution and we prove by induction on $u \geq a + b - 2$ that the $T_{a,b}(u)$ correspond
to the $q,t$-characters and that the powers of $t$ are given by the distinguished choice.
As discussed above, the $X_{a,b}$ are algebraically independent so we can identify the $T_{a,b}(u)$ 
for $u = a + b -  2$, with the corresponding $q,t$-characters. In general, we have a relation
$$T_{a,b}(U + 1) * T_{a,b}(U - 1) = t^\alpha  T_{a - 1,b}(U) * T_{a+1,b}(U) + t^\beta T_{a,b + 1}(U) * T_{a,b - 1}(U),$$
for some $\alpha, \beta\in\ZZ/2$. For $u\leq U$, we have $T_{a,b}(u) = M_{a,b}(u) \chi_{a,b}(u)$ where $M_{a,b}(u)$ is a monomial in the quantum torus and $\chi_{a,b}(u)$
is a polynomial in the $A_{i,c}^{-1}$ with coefficients in $\ZZ[t^{\pm 1}]$ and with constant term $1$ (see (\ref{order})). Then $(\chi_{a,b}(u))^{-1}$ is a formal power series in the $A_{i,c}^{-1}$. Each term of the sum 
$$T_{a,b}(U+1) = t^\alpha  T_{a - 1,b}(U) * T_{a+1,b}(U) *  (T_{a,b}(U-1))^{-1}+ t^\beta T_{a,b+1}(U) * T_{a,b-1}(U) * (T_{a,b}(U-1))^{-1}$$ 
is a monomial multiplied by such a formal power series.  The highest monomial is 
$$t^\alpha M_{a - 1,b}(U) * M_{a+1,b}(U) *  (M_{a,b}(U-1))^{-1}$$ 
which only appears in the first term. As $T_{a,b}(U+1)$ is bar invariant, it imposes that $\alpha$ is the power of the distinguished choice. Then one may consider 
$$T_{a,b}(U+1) - t^\alpha  T_{a - 1,b}(U) * T_{a+1,b}(U) *  (T_{a,b}(U-1))^{-1}.$$ 
The same arguments identifies $\beta$ with the distinguished choice. 
Hence $T_{a,b}(U+1)$ satisfies the equation as the corresponding $q,t$-character and so is equal to it.

\begin{ex}
Let us study the examples \ref{11}, \ref{12}, \ref{21} above. Let $n = \ell = 1$. We get : 
\begin{align*}L(Y_3) * L(Y_1) = t L(Y_1Y_3) + 1,\quad L(Y_1) * L(Y_3) = t^{-1}L(Y_1Y_3) + 1.\end{align*}
Let $n = 1$, $\ell = 2$. 
\begin{align*}
L(Y_3Y_5) * L(Y_1) = L(Y_1Y_3Y_5) + t^{-1} L(Y_5)
&, \quad L(Y_5) * L(Y_1Y_3) = L(Y_1Y_3Y_5) + t^{-1}L(Y_1),
\\L(Y_1) * L(Y_3) = 1 + t^{-1}L(Y_1Y_3)&, \quad  L(Y_1Y_3) * L(Y_3Y_5) = 1 + t^{-1} L(Y_1Y_3Y_5) * L(Y_3),
\\L(Y_3) * L(Y_5) = 1 + t^{-1}L(Y_3Y_5)&, \quad L(Y_3Y_5) * L(Y_1) = L(Y_1Y_3Y_5) + t^{-1} L(Y_5),
\\L(Y_5) * L(Y_1Y_3) = L(Y_1Y_3Y_5) + t^{-1}L(Y_1)&, \quad L(Y_1) * L(Y_3) = 1 + t^{-1}L(Y_1Y_3),
\\L(Y_1Y_3) * L(Y_3Y_5) = 1 + t^{-1}L(Y_3)*L(Y_1Y_3Y_5)&, \quad L(Y_3) * L(Y_5) = 1 + t^{-1}L(Y_3Y_5).
\end{align*}
Let $n = 2$, $\ell = 1$. 
\begin{align*}
L(Y_{1,3}) * L(Y_{1,1}Y_{2,4}) &= t^{\frac{1}{2}}L(Y_{2,4}) * L(Y_{1,1}Y_{1,3}) + L(Y_{2,2}Y_{2,4}),\\
L(Y_{2,4}) * L(Y_{1,1}) &= t L(Y_{1,1}Y_{2,4}) + 1,
\\L(Y_{1,1}Y_{2,4}) * L(Y_{2,2}) &= t^{\frac{3}{2}} L(Y_{1,1}) * L(Y_{2,2}Y_{2,4}) +  L(Y_{1,1}Y_{1,3}),
\\L(Y_{1,1}) * L(Y_{1,3}) &= t^{\frac{1}{2}} L(Y_{2,2}) + t^{-\frac{1}{2}}L(Y_{1,1}Y_{1,3}),
\\L(Y_{2,2}) * L(Y_{2,4}) &= t^{\frac{1}{2}}L(Y_{1,3}) + t^{-\frac{1}{2}}L(Y_{2,2}Y_{2,4}), 
\\L(Y_{1,3}) * L(Y_{1,1}Y_{2,4}) &= t^{\frac{3}{2}}L(Y_{1,1}Y_{1,3}) * L(Y_{2,4})  + L(Y_{2,2}Y_{2,4}),
\\L(Y_{2,4}) * L(Y_{1,1}) &= t L(Y_{1,1}Y_{2,4}) + 1, 
\\L(Y_{1,1}Y_{2,4}) * L(Y_{2,2}) &= t^{\frac{1}{2}} L(Y_{2,2}Y_{2,4}) * L(Y_{1,1})  +  L(Y_{1,1}Y_{1,3}),
\\L(Y_{1,1}) * L(Y_{1,3}) &= t^{\frac{1}{2}} L(Y_{2,2}) + t^{-\frac{1}{2}} L(Y_{1,1}Y_{1,3}), 
\\L(Y_{2,2}) * L(Y_{2,4}) &= t^{\frac{1}{2}}L(Y_{1,3}) + t^{-\frac{1}{2}}L(Y_{2,2}Y_{2,4}).
\end{align*}
\end{ex}

\end{document}